# Path Optimization and Near-Greedy Analysis for Graph Partitioning: An Empirical Study

Jonathan Berry[*] Mark Goldberg[†]

## Abstract

This paper presents the results of an experimental study of graph partitioning. We describe a new heuristic technique, *path optimization*, and its application to two variations of graph partitioning: the *max_cut* problem and the *min_quotient_cut* problem. We present the results of computational comparisons between this technique and the Kernighan-Lin algorithm, the simulated annealing algorithm, the *FLOW*-algorithm of [17], the multilevel algorithm of [14], and the recent 0.878-approximation algorithm of [7]. The experiments were conducted on two classes of graphs that have become standard for such tests: random and random geometric. They show that for both classes of inputs and both variations of the problem, the new heuristic is competitive with the other algorithms, and holds a advantage for *min_quotient_cut* when applied to very large, sparse geometric graphs (10,000 - 100,000 vertices, average degree $\leq 10$).

In the last part of the paper, we describe an approach to analyzing graph partitioning algorithms from the statistical point of view. Every partitioning of a graph is viewed as a result achieved by a "near greedy" partitioning algorithm. The experiments show that for "good" partitionings, the number of non-greedy steps needed to obtain them is quite small; moreover, it is "statistically" smaller for better partitionings. This led us to conjecture that there exists an "optimal" distribution of the non-greedy steps that characterize the classes of graphs that we studied.

## 1 Introduction

Given a graph $G = (V, E)$ and a partitioning $\pi$ of $V$ into disjoint sets $S$ and $\bar{S}$, let $C(\pi)$ denote the number of edges[1] cut. The goal of the *max_cut* problem is to find a partitioning $\pi$ which maximizes $C(\pi)$. The *quotient cost* of $\pi$ is defined as $\frac{C(\pi)}{min(|S|,|\bar{S}|)}$. Finding a cut with minimum quotient cost is called the *min_quotient_cut* problem. Both problems are known to be NP-hard. The problems have received a great deal of attention because of their applications, most notably in VLSI design (see [1]), and their potential usefulness for many other optimization problems ([3], [15], [7], [4], [20], [16], [17], [10], [5]). An experimental study of a heuristic algorithm for *min_quotient_cut* based on the multicommodity flow technique was done in [17]; the best approximation algorithm for *max_cut*, one with a provable approximation ratio of .878, was recently described in [7].

In §2, §3, and §4, we describe a new heuristic technique and its application to *max_cut* and *min_quotient_cut*. We present empirical comparisons between the new algorithm and the Kernighan-Lin [2] algorithm [16] (*KL*), the simulated annealing algorithm of [15] (*SA*), the *FLOW*-algorithm described in [17], and the 0.878-approximation algorithm given in [7]. The experiments suggest that the new algorithm is competitive with those algorithms, and that it is superior to them for certain classes of inputs.

A description of an approach to analyzing graph partitioning algorithms from a statistical point of view is contained in §5. A partitioning of the vertex set of a given graph is viewed as the result of a process which successively places vertices of the preordered vertex set into the "left" or "right" partitions. Each placement is labeled "greedy" or "non-greedy," depending on the number of additional edges cut. Given class $\mathcal{T}$ of graphs, we consider a function $F_{\mathcal{T}}(i)$, called the *ng*-function of the class, defined to be the probability that for a graph $G \in \mathcal{T}$, the $i^{th}$ vertex placement is non-greedy. Extensive experiments approximating $F_{\mathcal{T}}(i)$ for several classes of graphs are presented in §5. It turns out that for all the classes of input graphs considered, the partitionings constructed by the best heuristics contain a surprisingly small portion of non-greedy steps, and most of these are located among the "first" placements. Furthermore, we discovered that for a given input graph, better partitionings contain fewer non-greedy steps. This leads us to the conjecture that there exists an "optimal" *ng*-function corresponding to an optimal partitioning.

In the case of random graphs with $n$ vertices and edge probability $p$, the *ng*-function is approximated by the expression $a(n,p) + b(n,p)/\sqrt{i}$ $(i = 1, \ldots, n)$. Linear regression analysis indicates that $a$ and $b$ are nearly constant functions of $n$ and $p$. In turn, we use these statistics to derive a ***probabilistic greedy*** algorithm, *pg*-procedure, which produces the output based solely on the statistics. The experiments show that the partitionings produced by such a simplified procedure

[*]Resselaer Polytechnic Institute

[†]Rensselaer Polytechnic Institute; the work of this author was supported in part by NSF Grant #CCR-9214487.

[1]The technique described in this paper can be easily expanded to graphs with weighted edges.

[2]Our implementation uses the version of Fiduccia and Mattheyses (*FM*) [5]

are reasonably close to the best partitioning constructed by other algorithms.

## 2 Path Optimization

Many applications of partitioning concern *hypergraphs* rather than graphs. A hypergraph $G$ is a pair $(V, E)$, where $V = \{v_0, v_1, \ldots, v_{n-1}\}$ is a set of vertices and $E = \{e_0, e_2, \ldots, e_{m-1}\}$ is a collection of subsets of $V$, called[3] edges (or hyperedges). Although all of our experiments were conducted on graphs (hypergraphs with all edges being sets of size two) we describe $PO$ as a partitioning algorithm for hypergraphs. In the real code, a few simplifications can be made to take advantage of the fact that the inputs are graphs.

Starting with an initial partitioning, $PO$ constructs[4] a sequence of partitionings with non-decreasing values of the objective function; the algorithm halts if five partitionings do not produce an improvement in the objective function. Given a partitioning $(S, \bar{S})$, each iteration of $PO$ constructs a sequence $P$ of vertices that alternate between $S$ and $\bar{S}$. When the construction of the sequence is done, its vertices are "flipped-flopped," i.e., those vertices in $S$ are moved to $\bar{S}$ and vice-versa. The sequence $P$ is developed vertex by vertex. To make the construction efficient, we restrict the set of candidates for expanding $P$ to the vertices that are not in $P$ and are adjacent to the latest addition to $P$. The cost for flip-flopping $P$, which we call the *flip_cost* of $P$, is computed in increments vertex by vertex. The next path vertex is chosen to be the first candidate which does not make the *flip_cost* worse, i.e. smaller in the case of maximization and bigger in the case of minimization. The path development can be shown to take time linear in the number of edges, assuming the average degree and average edge size are bounded by constants.

Given a partitioning $(S, \bar{S})$, let $\text{loc}(v)$ denote the partition where $v$ resides. An edge $e$ containing $v$ is called type-0 critical (resp. type-1) with respect to $v$, if $\forall w \in e[(w \neq v) \Rightarrow (\text{loc}(v) = \text{loc}(w))]$ (resp. $\text{loc}(v) \neq \text{loc}(w)$). Let $n(v)$ (resp. $\bar{n}(v)$) be the number of type-0 critical (resp. type-1) edges with respect to $v$. The *gain* of a vertex $v$ denoted $cg(v)$ is defined as $n(v) - \bar{n}(v)$.

Below, we give an intuitive explanation of the $PO$ algorithm. The idea is simple, but the details are tedious; they will be presented in the full paper. Suppose that some sequence of vertices $P$ has been selected such that the locations of $w \in P$ alternate between $S$ and $\bar{S}$. Furthermore, suppose that for all $w \in P$, $w$ is marked as "locked." To describe the method of selecting the next vertex in the sequence, we use a function *flip_cost_incr* which, given a new vertex $v$, determines the change in *flip_cost*($P$) if $v$ were to be added. Pseudocode for all functions described here will be provided in the full paper.

The *flip_cost_incr* algorithm identifies edges of certain classes which introduce differences between $\sum_{v \in P} \text{cg}(v)$ and *flip_cost*($P$), and computes the increment to the *flip_cost* obtained by concatenating $v$ to the sequence $P$. It returns *TRUE* if the increment is non-negative in the case of maximization and non-positive in the case of minimization. A simple running time analysis shows that *flip_cost_incr* runs in $O(ds)$, where $d$ denotes the maximum degree of the graph, and $s$ denotes the maximum size of any edge. For very large graphs, practicality demands that $d$ and $s$ be bounded by constants, at least on average.

Let $v_i$ denote the most recent addition to $P$, and $v_{i-1}$ denote the second most recent. If the problem is minimization, the next path vertex is selected from the neighbors of $v_{i-1}$ which share its partition. In the case of maximization, the selection is from the neighbors of $v_i$ which are located in the opposing partition. The selected vertex is the first such that *flip_cost_incr* is not unfavorable. The running time of this selection algorithm is thus $O((ds))$,

Note that the next path vertex could just as easily be selected in a greedy way using a routine which returns the vertex with best cg such that the path is extended. However, this has been implemented and found not to perform any better than the simply adjacency list traversal.

$PO$ employs a greedy randomized initial partitioning generator, to be described below. Thus, the algorithm can be run from an arbitrary number of randomized starting partitions with the expectation that the best solution obtained improves with time, like $KL$ and $SA$. Given a partitioning, the $PO$ driver algorithm repeatedly finds a path with beneficial *flip_cost* and flip-flops it. Each iteration of this process begins by examining the $k$ vertices with best $cg$, looking for one which starts an alternating path with satisfactory *flip_cost.* If no such path is found, the iteration terminates. After five iterations terminate without an improvement to the objective function, the best partitioning is stored and the next initial partitioning is obtained. The cost of updating cell gains for all vertices is shown to be linear in the number of edges in [5]. The time complexity of an iteration of $PO$ is thus $O(m)$, assuming that the degree and edge size are bounded by constants. However, each iteration is much faster in practice. The average path

[3] In the VLSI literature, the vertices are often called *cells* and edges are called *nets*.

[4] Many versions of the algorithms were considered and tested; the one described below performs best on the classes of inputs we used.

length is usually less than 3 for very sparse graphs.

The *PO* algorithm is more closely analogous to greedy local optimization than the *KL* family of algorithms. Instead of developing a path which eventually includes all vertices and choosing the best intermediate swap state, *PO* remains a greedy hill climbing algorithm. However, vertex moves which would be considered non-greedy by local optimization, causing the process to halt, can be accepted because the "greediness" of a move is no longer associated with a single vertex. It seems that the *adjacencies* between swapped vertices are as important or more important than the individual cell gains.

The initial partitioning for *PO* is generated by a constructive greedy algorithm which starts with empty partitions, selects vertices one by one, and places them into partitions in a greedy way with respect to the current objective function. The vertex selection algorithms, described below, are *max-diff* for maximization problems, and *min-diff* for minimization problems. For each vertex placement, let $U$ be the set of unplaced vertices. Let $\delta(v) = |n(v) - \bar{n}(v)|$, $M_{\min} = \{v : \forall w \in U, \delta(v) \leq \delta(w)\}$, and $M_{\max} = \{v : \forall w \in U, \delta(v) \geq \delta(w)\}$. In *max-diff* selection, next vertex is drawn at random from $M_{\max}$, while *min-diff* selects one at random from $M_{\min}$. A similar vertex reordering technique is used and analyzed in [4]. Thus, each process is randomized, and the initial partitionings of *PO* are generated by random walks down the implicit backtracking tree. We refer to this initial partitioning generator as the $W$ algorithm.

## 3 Algorithm Implementations and Setup of Experiments

Path Optimization was compared extensively with the Kernighan-Lin (*FM* version) algorithm and the simulated annealing algorithm as described in [15]. Additional comparisons were made with the *FLOW* heuristic described in [17], the multilevel algorithm of [14], and the .878-approximation algorithm for *max_cut* [7].

The types of inputs studied were random graphs and random geometric graphs. The latter are created by laying out $n$ random points on the unit square and connecting only those whose Euclidean distance is less than a given threshold $d$. Both of these types of graphs have been used for comparisons before in the literature, and we continue the trend in order to facilitate further comparisons. Geometric graphs present to partitioning algorithms quite a different challenge from random graphs. In fact, the ranking of algorithms can be reversed when moving from random to geometric graphs, as we will see below.

The *KL*, *SA*, and *PO* algorithms are all randomized, and continue to perform iterations from different starting configurations as long as running time permits. The *PO* initial partitioning generation is described in §2. The *KL* and *SA* algorithms start from random partitionings, except when applied to the *min_quotient_cut* problem on geometric graphs. Here, initial partitionings are generated by the *line* heuristic described in [15]. The *line* heuristic uses geometric information to split the vertex set of a geometric graph into two equal sized halves with a line of randomly chosen slope. It has been demonstrated that such initial partitionings dramatically improve the performance of *KL* and *SA* [15, 17].

Our implementation of *SA* follows that of Johnson, et al. with one important exception: for *max_cut* on random graphs, better results are obtained if the running time is spread over one long annealing run instead of several shorter ones. However, for *max_cut* on geometric graphs and *min_quotient_cut* on both graph types, the cooling ratio is set as in [15] and iterations are performed until the time is up.

Our implementations of the algorithms support various objective functions, including those of *max_cut* and *min_quotient_cut*. The modifications to achieve this are small. Our version of *KL* was tested on the set of geometric graphs from [17], and it reported results comparable to those of their *KL* implementation, which in turn had been tested against that of [15]. *SA* was not tested against any previous data sets, but our implementation is based directly on [15] and reports similar results when run on similar inputs.

For each algorithm, we computed the running time spent in the main loop only. The input and initialization times were not included. This gave a slight advantage to *SA*, which first makes a trial run to correctly set its cooling ratio variable. The time was taken with the Unix *getrusage()* command. According to our experiments, the amount of work done per given time is virtually independent of the system load.

All trials involving graphs of less than 100,000 vertices were run on Sparc 10 machines with 44.2 SPECint92 ratings; the trials involving graphs of 100,000 vertices or more were run either on Sparc 10 machines rated at 65.2 SPECint92, or RS6000 machines rated at 117 SPECint92 (comparisons are only drawn between runs on the same type of machine).

For each variation of graph parameters, a data set of more than thirty graphs was generated if the number of vertices was less than 100,000, and each algorithm was run on all instances. Graphs with larger numbers of vertices were grouped into samples of size ten. For each algorithm, only the best solution for each graph was retained. The sample mean and standard deviation of this set of observations were then computed, as well as a 99% confidence interval for the true mean solution. For a discussion of $100(1-\alpha)\%$ confidence intervals, see [2]. In standard statistical practice, a confidence interval derived from a sample size of more than thirty trials allows an appeal to The Central Limit Theorem and an argument that, with a given confidence, the true mean lies somewhere in the interval, regardless of the distribution of the individual trials. If the number of trials is less than thirty, as with our experiments with graphs of 100,000 vertices or more, the confidence interval is obtained using the Students $T$ distribution and the assumption is made that the population of individual trials is normally distributed.

## 4 Results of Comparisons

Here we present results for relatively dense (aver. degree 110) random graphs, and relatively sparse (aver. degree $\approx$ 10) geometric graphs. A brief summary of the experimental results on large, sparse geometric graphs follows:

- The multilevel spectral algorithm of [14] is significantly faster than the others. Given equal time, it comes to dominate them as graph size increases (see Table 4). However, it is not known whether, given increasing time, its partitioningss can match those of *W-PO* and the other heuristics.
- *W-PO* has an advantage over the others when available running time is moderate. The advantage increases as graph density decreases. (see Figures 5 and 6). As graphs become sparser, the constructive greedy *W* algorithm provides better initial partitionings than the *line* heuristic. Combinations of *W-KL* and *W-SA* have not yet been examined, but will be in the full paper.
- *FLOW* has found the best quotient cut for a single, large graph. However, since it is a global method, no partitioning is produced until the algorithm completes. The running time requirements are much larger than those of the others.

The 99% confidence intervals for solution quality are presented in Figures 1-6 for *KL* ($\mu_\kappa$), *SA* ($\mu_\sigma$), and *PO* ($\mu_\rho$). Note that the scales of the graphs in the figures are adjusted to highlight the differences in the performance of the heuristics. The differences are small in absolute terms. This leads to very small percentage differences for the *max_cut* sample means. However, the solutions to *min_quotient_cut* have few enough cuts that *W-PO* can hold an advantage of 20% or more over *line-KL* and *line-SA* by cutting slightly fewer edges (see Table 3 in §4.2).

In Figures 1 and 2, the horizontal axis represents the cut percentage ($\frac{C(\pi)}{|E(G)|}\times 100$). Figures 3, 4, 5, and 6 plot the quotient cut ($\frac{C(\pi)}{\min(|S|,|\bar{S}|)}$). In all of these figures, the vertical axis has no significance. Table 3 presents the same results as Figure 6, showing the number of cut edges instead of the quotient cost.

**4.1 Max Cut** The results of Johnson et al. [15], which concerned Graph Bisection, a minimization problem, suggested that simulated annealing was slightly better than *KL* for random graphs and clearly worse for geometric graphs. Our results show that this is not the case in general, even for minimization of the cut for geometric graphs (see §4.2). In fact, for the case of *max_cut*, *SA* was the overall winner for both random and geometric graphs.

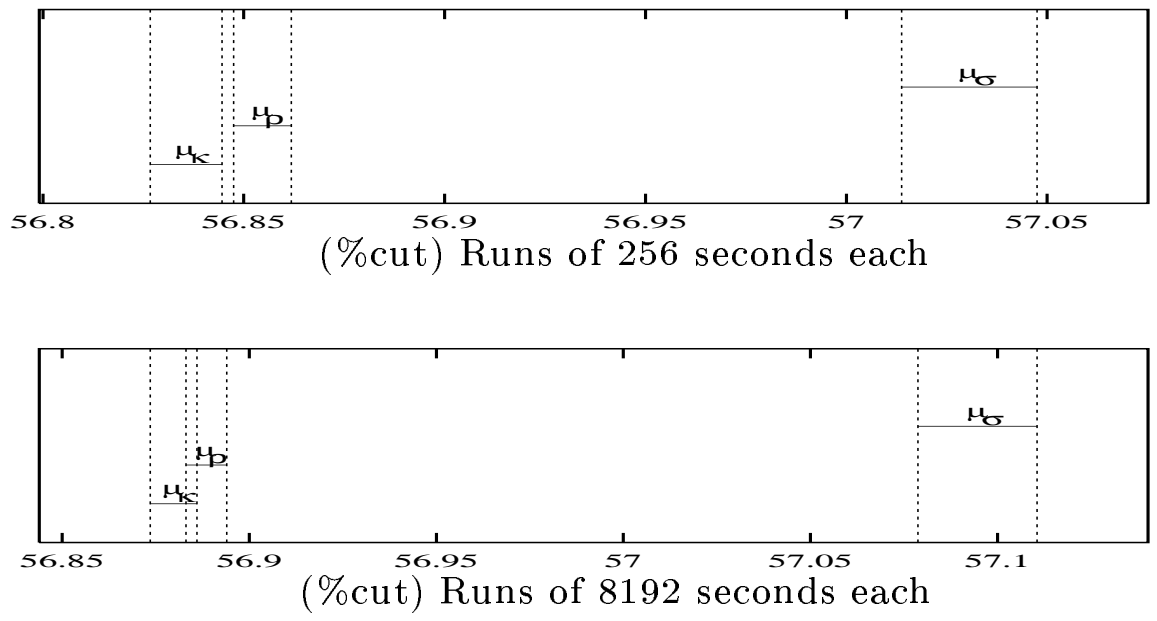

Figure 1: *max_cut*: 99% Confidence Intervals, Random Graphs, $n$ : 10,000, ave deg: 110

The results of our comparisons of *PO*, *KL*, and *SA* are presented in Figures 1 and 2. Although the confidence intervals are in most cases sufficient to rank the algorithms, with *SA* in the lead, the percentage difference between the average best solutions is so small that the separation between the various algorithms is not very significant.

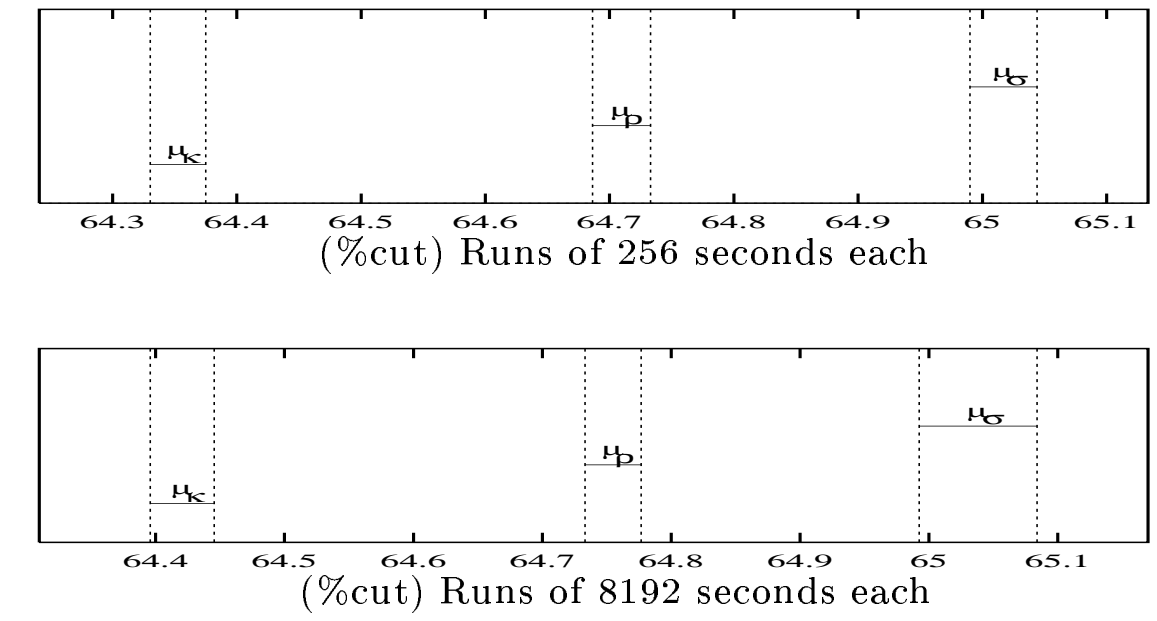

Figure 2: *max_cut*: 99% Confidence Intervals, Geometric Graphs, $n$ : 10,000, ave deg: 10

Recently, Goemans and Williamson ([7],[8]) constructed an approximation algorithm which delivers so-lutions to *max_cut* with a performance expectation of at least .87856 and also computes an upper bound which does not exceed the optimal value by more than a factor of $\frac{1}{0.87856}$. In Table 1, we compare the performance of *KL*, *SA*, *PO*, and their algorithm, *GW*, on two random graphs with 500 vertices and edge probabilities 0.05 and 0.5 respectively, and two random geometric graphs with 500 vertices each and distance thresholds of 0.05 and 0.5 respectively. Table 1 contains the ratios $c/u$, where $c$ is the value of the cut achieved by the corresponding algorithm, and $u$ is the upper bound computed by the program implementing the Goemans and Williamson algorithm. Note that obtaining upper bounds for larger sparse graphs is at the moment computationally infeasible. See the journal version of Geomans and Williamson's paper for details [8].

| | Random$(n,p)$ | | Geom$(n,d)$ |
|---|---|---|---|
| | (500,.05) | (500,. 5) | (500,.5) |
| *GW* | .9294 | .9783 | .9780 |
| *KL* | .9528 | .9844 | .9801 |
| *PO* | .9528 | .9842 | .9801 |
| *SA* | .9540 | .9852 | .9801 |

Table 1: *max_cut* bounds, 3 specific graphs

### 4.2 Min Quotient Cut

This section presents the results of our experiments with the *min_quotient_cut* problem. The modification needed to switch *KL* and *PO* to solve *min_quotient_cut* are straightforward. For *SA*, the balancing is achieved through a penalty function as in [15]. Our experiments with annealing based directly on changes in quotient cut offer no improvement in solution quality. When applied to *min_quotient_cut* on geometric graphs, our implementations of *KL* and *SA* employ the *line* heuristic described in §2. The initial partitioning of *PO* is obtained by the constructive greedy *W* algorithm also described in §2. It turns out that even though this greedy algorithm itself is beaten by *KL* and *SA*, it is an excellent initial partitioning generator, which is very effectively improved by *PO*. As graphs become very sparse (average degree $< 10$), the *W* algorithm produces better starts than the *line* heuristic.

Unlike the case of *max_cut*, there is a marked difference in the rankings of the algorithms between the random and geometric testbeds. The results presented in Figure 3 indicate that *SA* is still the best algorithm for partitioning these denser random graphs. However, for geometric graphs the rankings are reversed. Figure 4 shows the average best quotient cuts for the three heuristics for a test suite of fifty graphs of 10,000 vertices and average degree 10. In terms of the number of cuts achieved by the algorithms, we can be 99% confident for this graph generation process, size, density, and running time, and machine, the true expected number of cuts produced by our implementation of *W-PO* is between 131.5208 and 142.2792, while that of our implementation of *line-KL* is between 152.8726 and 163.5674. Thus, up to our assumptions, we can say with 99% confidence that for this case, our implementation of *W-PO* outperforms our implementation of *line-KL* by at least 6.9%. When the running time increases to 8192 seconds, the performance advantage increases to %10.3. In fact, after 50 runs of 8192 seconds each, the improved $\mu_\sigma$ and $\mu_\kappa$ do not match the sample mean of the 256 second runs of *W-PO*.

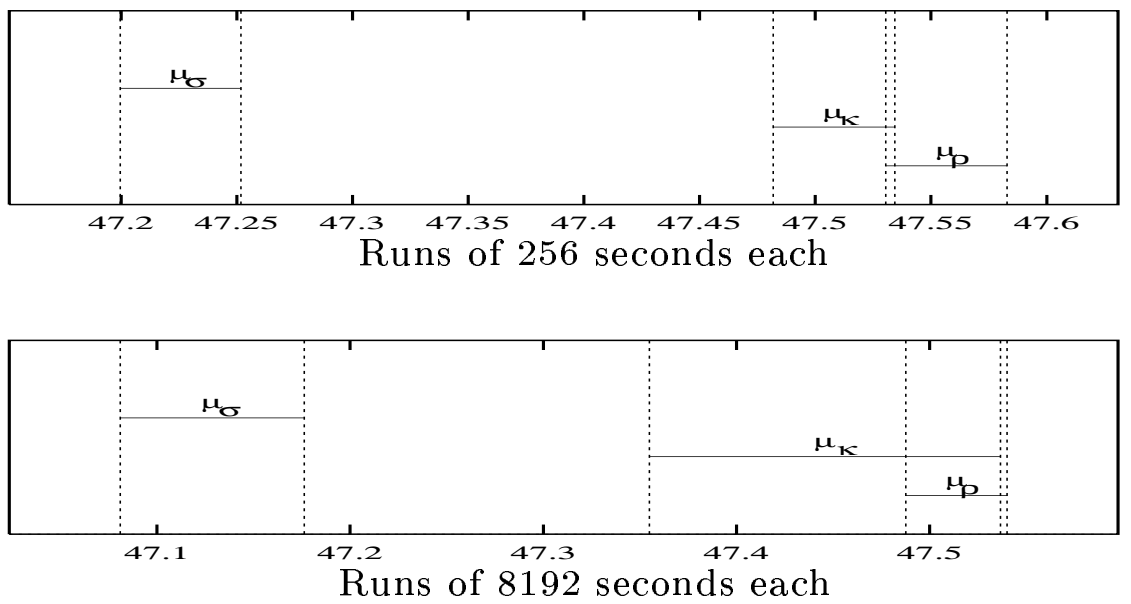

Figure 3: *min_quotient_cut*, Random Graphs, $n$ : 10,000, ave deg: 110

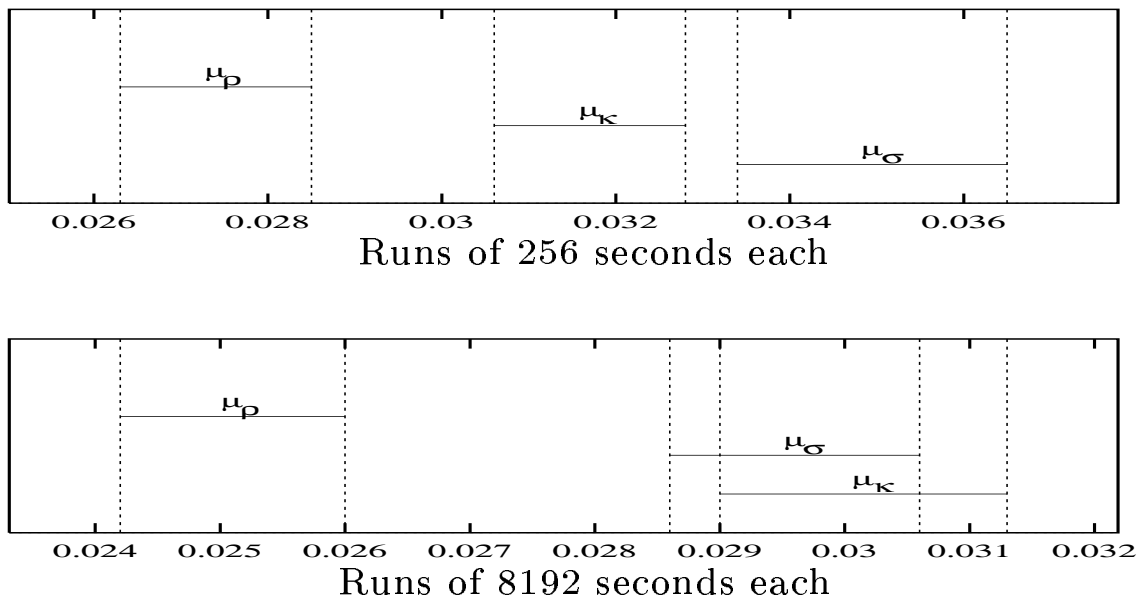

Figure 4: *min_quotient_cut*, Geometric Graphs, $n$ : 10,000, ave deg: 10

The question of scaling by graph size is addressed for geometric graphs of average degree 13.7 in Figure 5, and for those of average degree 7.6 in Figure 6 and Table 3. For the denser graphs of Figure 5, the running time of the algorithms is equalized at 10 hours for each graph. Again, all results which are presented in a given plot are taken from runs on the same type of machine. Since $\mu_\kappa$ has improved relative to $\mu_\rho$ such that the confidence intervals overlap, we cannot say with

| Size | Ave. Deg. | #graphs | line-*KL* | *FLOW-KL* | *W-PO* | #times *W-PO* best |
|---|---|---|---|---|---|---|
| 1000 | 10.94 | 10 | .1175 | .1134 | .1154 | 1 |
| 3000 | 11.46 | 5 | .0796 | .0726 | .0750 | 2 |
| 10000 | 12.02 | 5 | .0448 | .0406 | .0412 | 2 |

Table 2: *min_quotient_cut* Testbed of [17]

99% confidence that the true mean solution of either algorithm is better. However, the *sample* mean $\mu_\rho$ is still better than $\mu_\kappa$, and more trials tend to narrow the intervals.

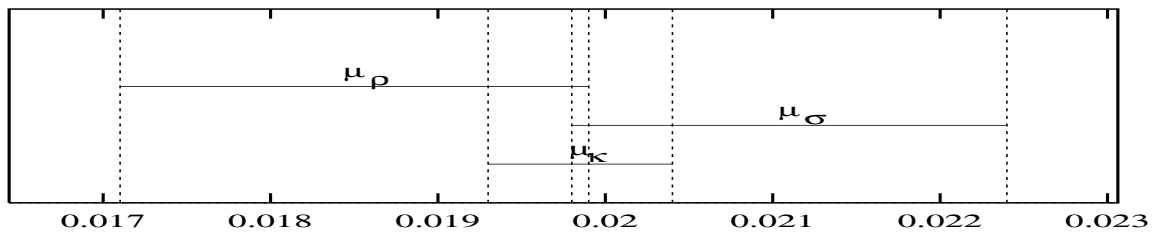

Figure 5: *min_quotient_cut*, 10 Geometric Graphs, $n$ : 100,000, Ave Deg: 13.7, Runs of 10 Hours

In [17], Lang and Rao described a heuristic, called *FLOW*, which uses the multicommodity flow approach to partitioning (see [18], [20]). They present the results of the empirical comparison of *FLOW* with variations of *KL* (*FM*) applied to sparse random and geometric random graphs. The authors conclud that *FLOW* is not useful for random graphs, but for geometric graphs, it achieves better results than *line-KL* as graph size increases, provided *FLOW* is augmented with *KL* or it is given longer running time. Our comparison of *FLOW* with *W-PO* for graphs of at most 10,000 vertices is given in Table 2. The *FLOW-KL* column refers to the quotient cut found by first applying *FLOW*, then cleaning up the solution with *KL*. For graphs of this size, there seems to be no clear winner. Although the average quotient cuts of *FLOW-KL* are slightly better, *PO* produced the best quotient cut five out of twenty times. In Figure 5, we considered a testbed of 10 geometric graphs of 100,000 vertices with average degree 13.7. In [17], *FLOW-KL* was run on one such graph. After 3 days of running on a 36 MHz silicon graphics machine, it produced a quotient cut of .014 ($\approx$ 700 cuts). In a run of similar duration on the same graph, the best achievable by *PO* was .019 ($\approx$ 950 cuts), and the best by *line-KL* was .020 ($\approx$ 1000 cuts). However, if running time is restricted to one hour, the cuts produced by *W-PO* and *line-KL* are only slightly worse.

Results for *min_quotient_cut* runs on a test suite of sets of sparser (ave. degree 7.6) geometric graphs of 12,500, 25,000, 50,000, 100,000, and 200,000 vertices are presented from two different points of view in Table 3 and Figure 6. These results, when compared to Figure 5, illustrate that *W-PO*'s advantage over *line-KL* and *line-SA* increases as graph density decreases, but decreases as graph size increases.

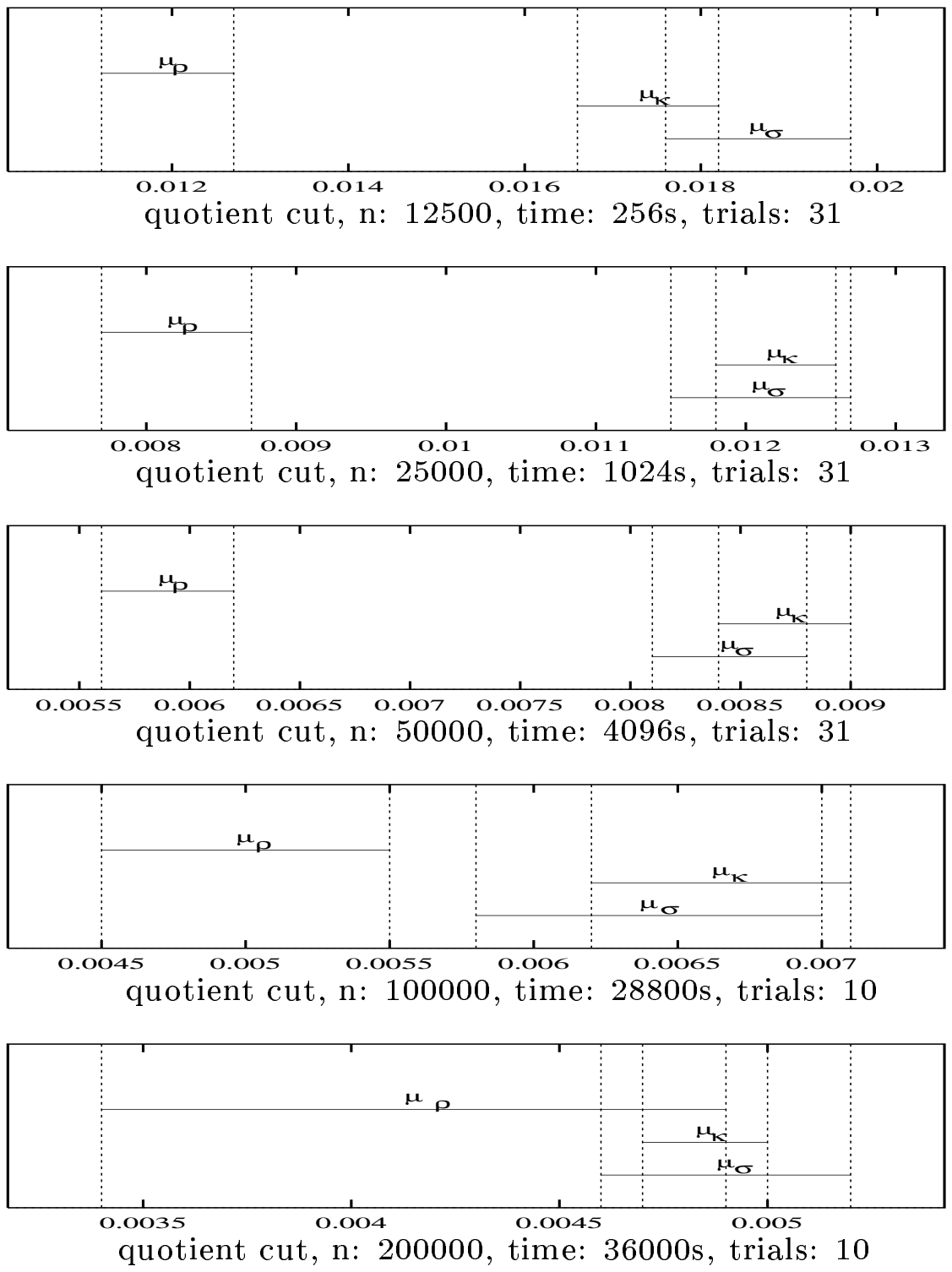

Figure 6: *min_quotient_cut*: 99% Confidence Intervals for the best quotient cut of *W-PO* ($\mu_\rho$), *line-KL*($\mu_\kappa$), and *line-SA* ($\mu_\sigma$), Geometric Graphs, Average Graph Degree: 7.6

Figure 6 shows the 99% confidence intervals for the expected quotient cost returned by our implementations of *W-PO*, *line-KL*, and *line-SA* for these sparse graphs. Table 3 gives confidence intervals for the expected number of cuts in the *min_quotient_cut* partitionings of the three algorithms. The expected advantage of *W-PO* over the others may appear relatively insignificant when quotient costs are examined; however, the expected advantage in cuts for this class of graphs is more than 20% in some cases. This advantage is determined by measuring the smallest interval between respective confidence intervals, and appears to decrease as graph size increases. However, even for the biggest graphs in the test set, the sample mean of *PO* holds a distinct advantage over the others. More trials would likely narrow the confidence intervals. An unanswered question is whether *W-KL* and *W-SA* can compete with *W-PO.*

| vertices | graphs | *W-PO* | *line-KL* | *line-SA* | *W-PO* min. exp. advantage |
|---|---|---|---|---|---|
| 12,500 | 31 | 70.51, 79.68 | 103.03, 113.16 | 110.45, 123.49 | 22.6% |
| 25,000 | 31 | 95.76, 109.02 | 147.81, 157.80 | 143.94, 159.09 | 24.2% |
| 50,000 | 31 | 139.90, 154.81 | 210.33, 224.38 | 203.53, 218.34 | 23.9% |
| 100,000 | 10 | 224.28, 273.52 | 310.49, 353.11 | 292.08, 349.32 | 6.3% |
| 200,000 | 10 | 337.75, 492.25 | 468.18, 500.02 | 464.91, 518.29 | no exp. advantage |

Table 3: 99% Confidence Intervals for Expected # Cuts, Geometric Graphs, average degree 7.6 (data from Figure 6)

| vertices | time | *W-PO* | *line-KL* | *line-SA* | *ML* |
|---|---|---|---|---|---|
| 12,500 | 4sec | 104.97,133.16 | 139.14,155.77 | 139.08,159.70 | 100.96,122.98 |
| 25,000 | 8sec | 158.79,193.27 | 194.31,218.60 | 190.59,218.25 | 158.37,194.08 |
| 50,000 | 17sec | 282.49,354.48 | 275.52,308.74 | 275.40,311.57 | 184.54,214.49 |
| 100,000 | 39sec | 350.83,869.37 | 388.12,548.68 | 390.38,508.82 | 285.01,345.79 |
| 200,000 | 90sec | 961.58,1572.62 | 558.25,679.75 | 570.13,692.47 | NA |

Table 4: Very short runs: (*ML* is the multilevel algorithm of[14]) Expected # Cuts (data from Figure 6)

The idea of partitioning large graphs by performing a series of graph contractions has been explored in [9], [3], [15], and [14]. In [3] and [15], empirical evidence is presented indicating that a contracted version of *KL* can improve the algorithm both in speed and quality if the input graphs are very sparse. Bui [3] finds an advantage for regular graphs of a special class only if the degree is four or less. In [14], Hendrickson and Leland give a multilevel algorithm which uses weighted intermediate graphs to preserve good partitionings as the graph is uncontracted. The contractions are obtained by finding maximal matchings and identifying endpoints of matching edges. After the resultant graph is partitioned using the spectral method of [12], the original graph is restored through a series of uncontractions, with *KL* (*FM*) occasionally cleaning the partitioning. Results are presented indicating that for bisection of large, sparse graphs, this algorithm performs significantly better than spectral partitioning alone.

Using *Chaco* [13], a partitioning system due to Hendrickson and Leland which implements several spectral partitioning methods and the multilevel algorithm described above, we were able to make limited comparisons with *W-PO*, *line-KL*, and *line-SA*. The algorithms were run on the same set of graphs as in Table 3 and the results are presented in Table 4. The intended application for *Chaco* is the mapping of parallel computations, where speed is obviously extremely important. The multilevel algorithm is very fast, while the heuristics, especially *W-PO* require some time to work well. The advantage of the multilevel algorithm increases with graph size for these short running times. However, its expected solution quality falls short of the longer runs of *W-PO* (see Table 3). The most straightforward way to extend the running time of the multilevel algorithm is to grant the *FM* cleanup routine more time for randomized iterations. In future work, we hope to experiment with this possibility, and test the multilevel algorithm with longer runs and *PO* as a cleanup routine.

## 5 Near Greedy Analysis

In this section, we present empirical results showing some surprising statistical trends in the distribution of non-greedy steps needed to obtain near-optimal solutions to the *max_cut* problem. All experiments were conducted using *KL* on graphs with up to 1500 vertices. In order to accumulate representative statistics, it was necessary to apply *KL* to at least 1000 graphs per set of parameters. However, the experiments suggest that collecting statistics for much larger graphs might not be necessary.

### 5.1 Postprocessing

The purpose of the postprocessing is to determine how much the partitionings produced by the best known algorithms differ from those obtained by the greedy algorithm. Given a partitioning $(S, \bar{S})$ of a graph $G$, the postprocessing is accomplished by a constructive algorithm which places vertices into partitions determined by *KL*, and marks them "greedy," or "non-greedy." Our code descends from a postprocessing system developed by Goldberg and Fox [6]. Our experiments showed that the choice of the selection procedure has a dramatic effect on the data. When using no reordering (random reordering), approximately 20% of the placements are marked as non-greedy. For varied tests with the *max-diff/max-degree* reordering defined in §2, an average of about 4% placements were marked non-greedy.

We used five different types of randomly-generated graphs for our experiments: random, random geometric, random regular, $k_1, k_2$-unbalanced regular, and $p_1, p_2$-unbalanced random. The unbalanced graphs are formed by bisecting the vertex set, and assigning edges using different edge probabilities within the partitions and between them. For every of class of inputs and every combination of $n$ and $p$, at least 1000 randomly generated instances were run. For each $i = 0, \ldots, n-1$, the ratio $k/n$ was computed, where $k$ is the number of graphs in which the $i$th vertex placement was labeled "non-greedy;" $k/n$ is termed the non-greedy probability. In Figures 7 and 8, the distributions are displayed with the vertex placement *percentile* on the $x$ axis, and *ng*-function(i), the average non-greedy probability for vertices in percentile $i$, on the $y$ axis. For each box in the figures, the distributions for 500, 1000, and 1500 vertices have been plotted. They suggest the following conclusions:

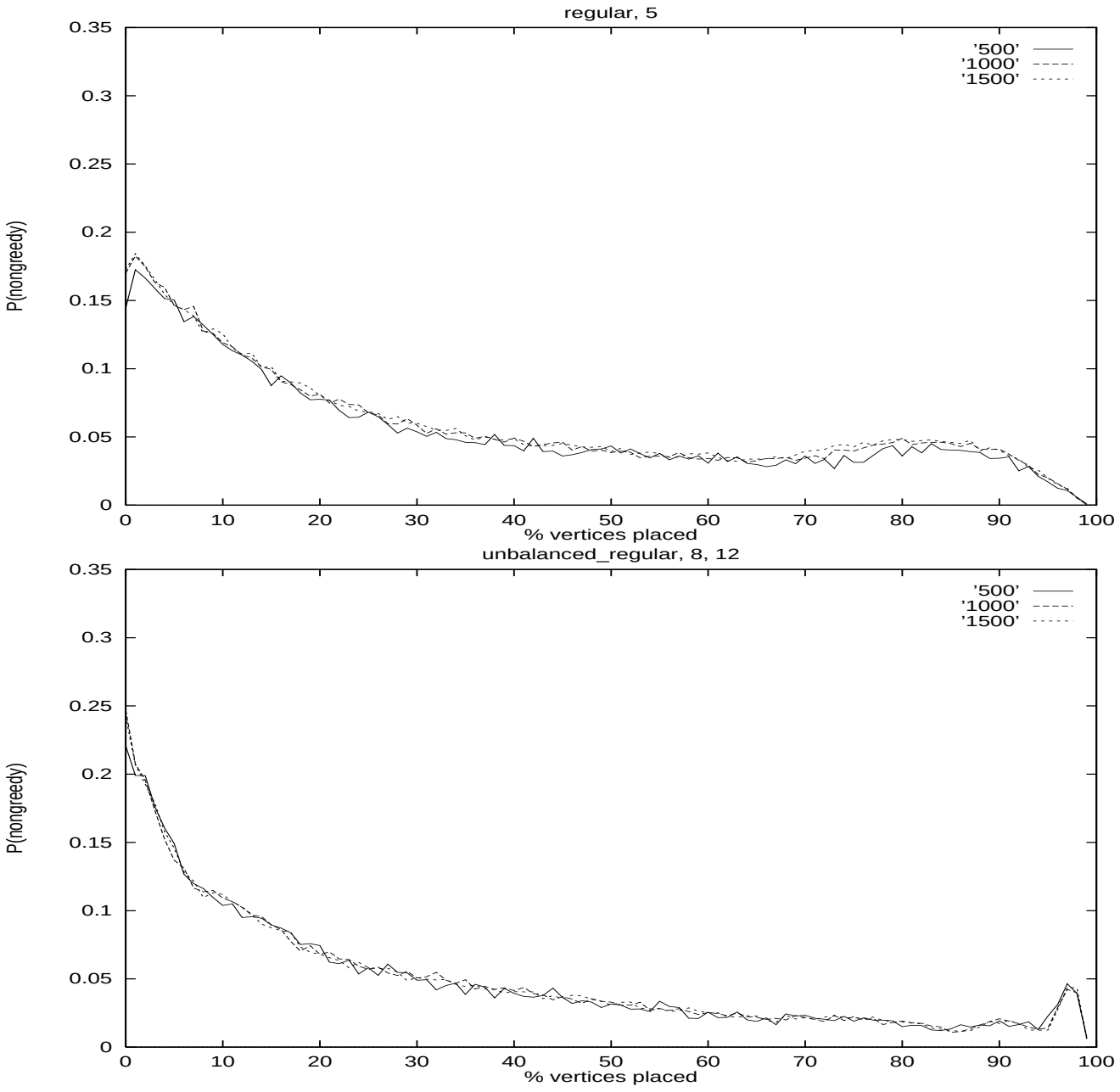

Figure 7: Nongreedy Distribution, Regular & Unbalanced Regular Graphs

1. For random and random regular graphs, the non-greedy distribution is virtually independent of the graph size.

2. There is a small constant $c$, such that for $i \geq c$, $ng(i)$ is almost a monotonically decreasing function for random and random regular graphs. For these graphs, approximately 80% of the non-greedy steps occur within the first 50% of the vertex placements. For geometric graphs, however, the non-greedy distribution is nearly a flat function.

3. The spiky increases in non-greedy probability for large arguments ( see Figures 7 and 8) occur when the average degree of the subgraph induced by the unplaced vertices is close to 1.

4. For random and regular graphs, the slope of the non-greedy distribution curve increases with graph density.

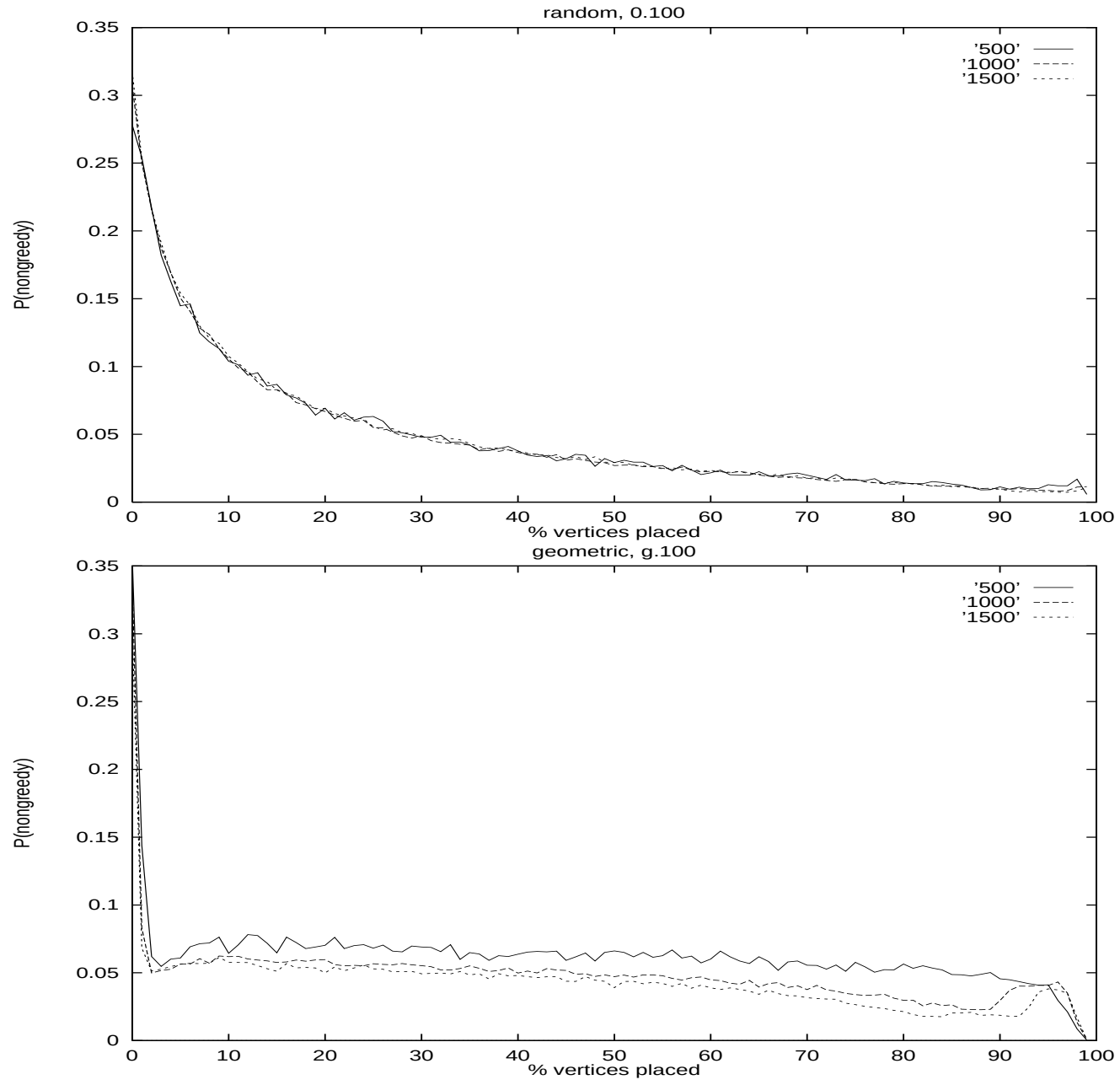

Figure 8: Nongreedy Distribution, Random and Geometric Graphs

Observation 2 can be explained intuitively by the graph type and the *max-diff/max-degree* reordering method used to construct the distributions.

The spikes noted in Observation 3 illustrate an increase in non-greedy probability that occurs near the end of the vertex placements in sufficiently sparse random and regular graphs. When the average degree of the subgraph induced by the unplaced vertices approaches 1.0, the non-greedy probability starts to increase, and it peaks when the average is slightly less than 1. We observed that the larger the average vertex degree of the graph, the later this occurs (see Figure 7 and in the top two boxes of Figure 8.)

For random graphs ($\mathcal{R}$(n,p), we used a linear regression model to approximate the *ng*-function $F_{\mathcal{R}}()$. The $R^2$ value shows the goodness of fit, i.e. how much of the variation (see [2]) of the dependent variable is explained by the independent variable. [5]

- $F_{\mathcal{R}_{500,.05}}(i) = -.0292 + \frac{.3991}{\sqrt(i)}$, $R^2 = .7645$, stderr $= .0490$
- $F_{\mathcal{R}_{500,0.5}}(i) = -.0369 + \frac{.4540}{\sqrt(i)}$, $R^2 = .9542$, stderr $= .0253$
- $F_{\mathcal{R}_{1900,.5}}(i) = -.0326 + \frac{.4294}{\sqrt(i)}$, $R^2 = .9249$, stderr $= .0302$

As reflected in the increasing $R^2$ values, we obtain better approximations as either graph density or graph size increases. The coefficients are nearly constant functions of $n$ and $p$ in absolute terms, but show a small increase in slope with increases in size or vertex degree.

**5.2 Probabilistic-Greedy Heuristic** The experiments with postprocessing the partitionings constructed by *KL*, led us to experiment with a "probabilistic greedy", *PG*, algorithmic strategy (see also, [10] and [11]). Under this paradigm, a solution is constructed successively, and for every step, the algorithm decides probabilistically if the step must be greedy or non-greedy. The decision is based on the *ng*-function which is supposed to be developed based on the previous experiments with the inputs of the specific class. Thus, the strategy is mainly a function of a class; the individual properties of a particular input come to play only when the two possible placements are classified as greedy or non-greedy. Since the procedure is probabilistic, it is repeated a number of times depending on the total running time available.

According to our experiments, given a small time allowance, *PG* lags behind both *KL* and the greedy heuristics. However, as the running time allotted to each increases, the *PG* solution tends to surpass that of the greedy algorithm and approach that of *KL*. Preliminary experiments show that there is a correlation between some simple parameters of graphs and the distribution of the non-greedy steps in its near-optimal partitionings. However, our attempts to use these correlations so far have not yielded a substantial improvement in the quality of the *PG* heuristic.

[5] The details of the model will be presented in the full paper.

## 6 Acknowledgments

We would like to thank Michael Goemans and David Williamson for providing us with a C-code of their .878 approximation algorithm (based on the code of Robert Vanderbei[19]), Bruce Hendrickson and Robert Leland for the *Chaco* system, and Satish Rao for sharing his test graphs. We are grateful to David Johnson for useful comments on Simulated Annealing, and to Vance Faber of Los Alamos National Laboratory for his support.